**Proofs of certain properties of irrational roots**

S. A. Belbas

Mathematics Department, University of Alabama, Tuscaloosa, AL 35487-0350.

*e-mail*: sbelbas@gmail.com

Abstract: We give two elementary proofs, at a level understandable by students with only pre-calculus knowledge of Algebra, of the well known fact that an irreducible irrational n-th root of a positive rational number cannot be solution of a polynomial of degree less than n with rational coefficients. We also state and prove a few simple consequences.

1. Introduction.

This paper is a set of class notes, on certain topics related to irrational numbers, proved in a manner that can be understood by students with only pre-calculus knowledge of Algebra.

Many algebra courses include proofs of the irrationality of $\sqrt{2}$, $\sqrt{n}$ if n is a natural number that is not a perfect square, and perhaps also proofs of results like the impossibility of having an irrational cubic root as a solution of a quadratic equation with rational coefficients. A collection of impossibility theorems (not restricted to irrational radicals) can be found in [R].

A well known result is that an irrational irreducible n-th root of a positive rational number cannot be solution of a polynomial equation of degree less than n with rational coefficients. (The terminology is explained in section 2 below.) Two elementary and simple proofs of this result are provided in these notes. A few consequences are also stated and proved.

It is expected that presentation of this material to students at the appropriate level can lead to an initial understanding of some basic concepts in the theory of equations, particularly the notion of irreducibility.

## 2. Irreducible radicals and a theorem on polynomials.

We start with a definition of irreducible radicals.

Definition 2.1. If A is a positive rational number, and n is a natural number with $n \geq 2$, we say that $\sqrt[n]{A}$ is irreducible if it is not possible to find a natural number m smaller than n and a positive rational number B such that $\sqrt[n]{A} = \sqrt[m]{B}$. ///

It follows from this definition that, if $\sqrt[n]{A}$ is irreducible, then it is irrational, because if $\sqrt[n]{A}$ were rational, say equal to a positive rational number R, then $\sqrt[n]{A} = \sqrt[m]{B}$ would hold with $m = 1$ and $B = R$. Also, for $n = 2$, the irreducibility of $\sqrt[n]{A}$ is equivalent to the irrationality of $\sqrt[n]{A}$.

We have

Proposition 2.1. For $n \geq 3$, the statement that $\sqrt[n]{A}$ is irreducible is equivalent to the statement that, for all natural numbers k satisfying $2 \leq k \leq n-1$, the number $\sqrt[n]{A^k}$ is irrational.

Proof: Suppose $\sqrt[n]{A}$ is irreducible. If, for some k with $2 \leq k \leq n-1$, we have $\sqrt[n]{A^k} = C$ with C a positive rational number, then $\sqrt[n]{A} = \sqrt[k]{C}$ which contradicts the irreducibility of $\sqrt[n]{A}$. The converse implication is proved in a similar way. ///

A characterization of irreducible radicals is given by the following

Proposition 2.2. If A and B are positive rational numbers, m and n are natural numbers greater than 1, and $\sqrt[n]{A} = \sqrt[m]{B}$ with $\sqrt[m]{B}$ irreducible, then m divides n and $A = B^{\frac{n}{m}}$.

Proof: Suppose $n = ms + r$ with s and r nonnegative integers and $1 \leq r \leq m-1$. The equality $\sqrt[n]{A} = \sqrt[m]{B}$ gives $A = B^s B^{\frac{r}{m}}$, so that $B^{\frac{r}{m}}$ is rational with r being a natural number less than m, which contradicts the irreducibility of $\sqrt[m]{B}$. ///

Of course, a partial converse is straightforward: if m divides n and $A = B^{\frac{n}{m}}$, then $\sqrt[n]{A} = \sqrt[m]{B}$.

It follows from Proposition 2.2 that the value $m_*$ that yields an irreducible representation of an irrational radical $\sqrt[n]{A}$ is the smallest among all positive divisors m of n for which $A^{\frac{m}{n}}$ is rational, and the value of B for which $\sqrt[n]{A} = \sqrt[m_*]{B}$ is $B = A^{\frac{m_*}{n}}$.

Now we shall prove:

Theorem 2.1. If $\sqrt[n]{A}$ is irreducible, then it cannot be a root of a polynomial p(x) of degree less than n with rational coefficients.

Proof: The idea for this proof is taken from the proof of a different result in [BP], article 238, where the idea is attributed to Abel.

Let $p(x) = x^m + a_{m-1}x^{m-1} + \cdots + a_1 x + a_0$ with m<n and rational coefficients $a_0, a_1, ..., a_{m-1}$. Assume that $\sqrt[n]{A}$ is a root of p(x). Set $q(x) := x^n - A$. Let d(x) be the greatest common divisor of p(x) and q(x). The standard Euclidean algorithm for finding the greatest common divisor of 2 polynomials (cf., e.g., [GN]) implies that d(x) has rational coefficients. The fact that, by assumption, p(x) and q(x) have at least one common root implies that the degree of d(x) is at least 1; the fact that d(x) exactly divides p(x) implies that the degree of d(x), say k, satisfies k<n. Suppose $d(x) = x^k + b_{k-1}x^{k-1} + \cdots + b_1 x + b_0$, with $1 \leq k<n$ and rational coefficients $b_0, b_1, ..., b_{k-1}$. Each of the roots of q(x) has modulus $\sqrt[n]{A}$, since all roots of q(x) are given by $r_\ell = \sqrt[n]{A} \exp\left(\frac{2\ell\pi i}{n}\right)$, $\ell = 0, 1, 2, ..., n-1$; since d(x) exactly divides q(x), the roots of d(x) form a subset of the roots of q(x); $|b_0|$ equals the product of the moduli of the roots of d(x), thus $|b_0| = \sqrt[n]{A^k}$, and by the irreducibility of $\sqrt[n]{A}$ and Proposition 2.1, $b_0$ is irrational, which contradicts the fact that all the coefficients of d(x) are rational. ///

Remark 2.1. The same contradiction obtains if p(x) is a polynomial with complex rational coefficients (i.e. complex numbers such that their real and imaginary parts are rational numbers), and we assume that p(x) and q(x) have at least one (complex) root in common. ///

3. Consequences of the previous results.

We state and prove a few direct consequences of the theorem of the previous section.

Consequence 3.1. If a, b, c, d, f, g are rational numbers, with c and g also positive, and if m and n are natural numbers greater than 1, if both $\sqrt[m]{c}$ and $\sqrt[n]{g}$ are irreducible, and if $a + b\sqrt[m]{c} = d + f\sqrt[n]{g}$, then $m = n, a = d$, and $b\sqrt[n]{c} = f\sqrt[n]{g}$.

Proof: We consider the case in which both b and f are non-zero (if at least one of b, f vanishes, the proof is quite simple). Without loss of generality, we assume that $m \leq n$. Then we obtain

$b^m c = \left((d-a) + f\sqrt[n]{g}\right)^m$. If m<n, then the expression $\left((d-a) + f\sqrt[n]{g}\right)^m - b^m c$ is a polynomial, with rational coefficients, of degree m<n in $\sqrt[n]{g}$, and therefore it cannot vanish according to theorem 2.1. Thus $m = n$. If $d \neq a$, by the binomial expansion of $\left((d-a) + f\sqrt[n]{g}\right)^n$, the equality $\left((d-a) + f\sqrt[n]{g}\right)^n - b^m c = 0$ implies that $\sqrt[n]{g}$ is a root of a polynomial of degree $n-1$, which is impossible by theorem 2.1. Thus $a = d$ and consequently $b\sqrt[n]{c} = f\sqrt[n]{g}$. ///

Consequence 3.2. If a, b, c, d are rational numbers, with b and d positive, if m and n are natural numbers greater than 1 with both $\sqrt[m]{b}$ and $\sqrt[n]{d}$ irrational, and if $a + \sqrt[m]{b} = c + \sqrt[n]{d}$, then $a = c$ and $\sqrt[m]{b} = \sqrt[n]{d}$.

Proof: Let b', d' be positive rationals and m', n' natural numbers greater than 1, such that $\sqrt[m']{b'} = \sqrt[m]{b}$ and $\sqrt[n']{d'} = \sqrt[n]{d}$, with both $\sqrt[m']{b'}$ and $\sqrt[n']{d'}$ irreducible. Then $a + \sqrt[m']{b'} = c + \sqrt[n']{d'}$ and by the preceding result we must have $a = c$ and $\sqrt[m']{b'} = \sqrt[n']{d'}$, consequently $\sqrt[m]{b} = \sqrt[n]{d}$.
///

Consequence 3.3. If A and B are positive rational numbers, m and n are natural numbers greater than 1, $\sqrt[m]{A}$ and $\sqrt[n]{B}$ are both irreducible, and either (i) $m \neq n$ or (ii) $m = n$ and $\sqrt[n]{\frac{A}{B}}$ is irrational, then the 3 numbers 1, $\sqrt[m]{A}$, and $\sqrt[n]{B}$ are linearly independent over the rational numbers.

Proof: It is to be shown that if a, b, c are rational numbers such that $a + b\sqrt[m]{A} + c\sqrt[n]{B} = 0$ then $a = b = c = 0$. If exactly one of b, c is zero, say $b = 0$ and $c \neq 0$, then we get $\sqrt[n]{B} = -\frac{a}{c}$ which contradicts the irrationality of $\sqrt[n]{B}$. If both b and c vanish, then $a = 0$. If both b and c are non-zero, then, by applying Consequence 3.1 to the equality $a + b\sqrt[m]{A} = -c\sqrt[n]{B}$ we obtain $a = 0$, $m = n$, and $b\sqrt[n]{A} = -c\sqrt[n]{B}$, thus $\sqrt[n]{\frac{A}{B}} = -\frac{c}{b}$ which contradicts the condition of irrationality of $\sqrt[n]{\frac{A}{B}}$. ///

(More general results, in the style of the above proposition, can be found in [B, S].)

Consequence 3.4. If A is positive rational and $\sqrt[n]{A}$ is an irreducible radical, then an expression of the form $a + b\sqrt[n]{A}$ with a and b rationals and $b \neq 0$ cannot be root of a polynomial of degree <n with rational coefficients.

Proof: If p(x) is such a polynomial, then the polynomial $p_1(x) := p(a + bx)$ is a polynomial with rational coefficients and has the same degree as p(x), thus it cannot have $\sqrt[n]{A}$ as a root, therefore p(x) cannot have $a + b\sqrt[n]{A}$ as a root. ///

4. Connections with irreducible polynomials.

In this section, "irreducible" (applied to polynomials) will mean irreducible over the rationals. Thus a polynomial p(x) with rational coefficients will be called *irreducible* if it cannot be factored as the product of two other non-constant polynomials with rational coefficients. Another bit of terminology is that a complex number is called *algebraic* if it is root of a polynomial with rational coefficients. A unique factorization theorem holds, with irreducible polynomials playing a role analogous to prime integers: each monic polynomial with rational coefficients can be uniquely factored as a product of monic irreducible polynomials. (A polynomial is termed *monic* if the coefficient of its highest degree term equals 1.) These properties are covered in practically every Algebra textbook, for example in [GN].

Two simple facts about irreducible polynomials are contained in the following two lemmas.

The first lemma extends the corresponding theorem of [D] about irreducibility of $x^p - A$ when p is a prime number and A is a positive rational number. On the other hand, more comprehensive results than lemma 4.1 below can be found in [L].

Lemma 4.1. If A is a positive rational number and $\sqrt[n]{A}$ is an irreducible radical, then the polynomial $x^n - A$ is irreducible.

Proof: If not, $x^n - A$ contains a monic polynomial factor $p_1(x)$ of degree say k, $1 \le k \le n-1$, with rational coefficients. The same argument as in the proof of theorem 2.1 shows that the zero-degree coefficient of $p_1(x)$ has absolute value $\sqrt[n]{A^k}$, which is impossible since the zero-degree coefficient of $p_1(x)$ is rational whereas $\sqrt[n]{A^k}$ is irrational by the irreducibility of $\sqrt[n]{A}$. ///

The next lemma is a well known result in the theory of algebraic equations; essentially equivalent statements can be found in [D] and many other textbooks. For the sake of completeness, we include the proof.

Lemma 4.2. If s is an algebraic number, then the monic irreducible polynomial with rational coefficients that has s as a root is unique.

Proof: Let p(x) and q(x) be monic irreducible polynomials with rational coefficients with $p(s) = q(s) = 0$, and let d(x) be the monic greatest common divisor of p(x) and q(x). The algorithm for finding a greatest common divisor implies that d(x) has rational coefficients. Since p and q have at least one root in common, the degree of d is at least 1. If the degree of d is less than the degree of p, then p can be factored as $p(x) \equiv d(x) p_1(x)$ where $p_1(x)$ has rational coefficients and the degree of $p_1$ is at least 1 and less than the degree of p(x), which is impossible when p(x) is irreducible. Thus the degree of d(x) equals the degree of p(x), therefore $p(x) \equiv d(x)$. Similarly $q(x) \equiv d(x)$. Consequently $p(x) \equiv q(x)$. ///

With the help of these two lemmas, we can give a second proof of theorem 2. This proof is very simple, given the two lemmas above, so simple in fact that some authors give the conclusion without proof, for example [B, P]. The theorem below is a particular case of the well known result (for instance, an exercise in [D]) that a root of an irreducible polynomial of degree n with rational coefficients cannot be a root of a polynomial of degree less than n with rational coefficients. Again, for completeness, we offer the few lines of proof.

<u>Second proof of Theorem 2.1:</u>  If p(x) is a monic polynomial with rational coefficients such that $p(\sqrt[n]{A}) = 0$, then p(x) contains a monic irreducible polynomial factor, with rational coefficients, which factor also has $\sqrt[n]{A}$ as a root. By the 2 lemmas above, that factor is $x^n - A$, which is impossible if the degree of p(x) is less than n. ///